\documentstyle[amstex,12pt]{article}

\def\m{\mbox{$\frak m$}}
\def\PP{\mbox{$\Bbb P$}}
\def\PN{\mbox{${\Bbb P}^N_K$}}

\def\QQ{\mbox{$\Bbb Q$}}

\def\codim{\mbox{\rm codim}}

\def\Proj{\mbox{\rm Proj}}

\def\reg{\mbox{\rm reg}}

\newtheorem{th}{Theorem}[section]

\newtheorem{prop}[th]{Proposition}
\newtheorem{lem}[th]{Lemma}

\newcommand{\qed}{$\Box$}

\begin{document}

\title
{Generic hyperplane section of curves
and \\ an application to regularity bounds \\
in positive characteristic}

\author{Edoardo Ballico \\
Department of Mathematics \\
University of Trento \\
38050 Povo (TN) \\
Italy \\
ballico@@science.unitn.it
\and
Chikashi Miyazaki \\
Nagano National College of Technology \\
716 Tokuma, Nagano 381-8550 \\
Japan \\
miyazaki@@cc.nagano-nct.ac.jp}
\date{}

\maketitle

\begin{abstract}

This paper investigates the Castelnuovo-Mumford regularity
of the generic hyperplane section of projective curves
in positive characteristic case, and yields an application
to a sharp bound on the regularity for nondegenerate
projective varieties.
\\

\noindent
{\it 1991 Mathematics Subject Classification.
Primary 14B15; Secondary 13D40, 14H45, 14N05} 
\end{abstract}

\section{Introduction}

The purpose of this paper is to study an upper bound of
the index of regularity of a generic hyperplane section
of projective curves and its
application to sharp regularity bounds for projective
varieties.
\\

For a projective scheme $X  \subset   \PN$,
we define the Castelnuovo-Mumford regularity $\reg (X)$
as the smallest integer $m$ such that
${\rm H}^i({\PN}, {\cal I}_X(m-i))  =  0$
for all $i  \ge  1$, see, e.g., \cite{E}.
The interest in this concept stems partly from the well-known fact:
The regularity $\reg (X)$ is the smallest integer $m$ such that
the minimal generators of the $n$-th syzygy module of the defining
ideal $I$ of $X$ occur in degree 
$\le  m  +   n$ for all $n \ge 0$. 
\par

In particular, for a zero-dimensional scheme $S  \subset   \PN$,
we define the index of regularity
$i(S)$  of  $S$  as the smallest integer $t$ such that
${\rm H}^1(\PN, {\cal I}_S(t))  =  0$.
We remark that $\reg (S)  =   i(S) +  1$.
\\

Throughout this paper, for a rational number
$\ell   \in    {\QQ}$, 
we write $\lceil  \ell  \rceil$
for the minimal integer which is larger than or equal to
$\ell$, and $\lfloor \ell \rfloor$
for the maximal integer which is smaller than or equal to
$\ell$.
\\

Let $S \subset \PN$ be a generic hyperplane section of
a nondegenerate projective curve $C \subset {\PP}_K^{N+1}$
over an algebraically closed field $K$.
Then $S$ has the uniform position property in case
${\rm char}(K) = 0$, see \cite{HE},
while the property does not necessarily hold in case
${\rm char}(K) > 0$, see \cite{R}.
Instead, even for the positive characteristic case,
$S$ has the linear semi-uniform
position property introduced in \cite{B1},
see $\S 2$ for the definition.
The linear semi-uniform position has an important
role in studying the positive characteristic case.
\par

For example, by studying the $h$-vectors
of a zero-dimensional scheme $S$ in linear semi-uniform
position, we have
an upper bound on the index of regularity, that is,
$i(S) \le \lceil
 (\deg (S) - 1)/N  \rceil$, see, e.g., \cite{B1,NS2}.
Also, there are some known facts on the sharpness of the above
bound.
If a zero-dimensional scheme
$S \subset \PN$ lies on a rational normal curve, then
we have an equality,
$i(S) = \lceil (\deg (S)  -  1)/N  \rceil$.
On the other hand, we assume that a zero-dimensional scheme
$S \subset \PN$ is in uniform position and
$\deg (S)$ is large enough.
If the equality 
$i(S) = \lceil (\deg (S)  -  1)/N  \rceil$ holds,
then $S$ lies on a rational normal curve, see, e.g., \cite{M,Y2}.
\\

In Section~2, we consider a generic hyperplane section
$S \subset \PN$ of
a nondegenerate projective curve over an algebraically closed
field $K$ such that $S$ does not have
the uniform position property. So we always focus on the case
${\rm char}(K) > 0$.
First, we will show that, under the condition that
 $N \ge 3$ and $\deg (S)$ is large
enough, if $S$ does not have the uniform position property,
then
$i(S) \le \lceil (\deg (S) - 1)/N \rceil - 1$ in (\ref{lem-v(1)=3})
and (\ref{lem-v(2)=4}).
The lemmas are technically key results of this paper.
As in classical Castelnuovo's method, we will show
the assertion of the lemmas, and in fact, the
linear semi-uniform position property will be useful
for this proof.
Then we apply the lemmas to the main result of this section, see Theorem~\ref{main-th}.
Let $S \subset \PN$ be a generic hyperplane section
of a nondegenerate projective curve with $\deg (S)$
large enough.
Without assuming $S$ is in uniform position, if the equality 
$i(S) = \lceil (\deg (S)  -  1)/N  \rceil$ holds,
then $S$ lies on a rational normal curve.
Finally we describe a results on the index of regularity
for a generic hyperplane section of very strange curves,
see Proposition~\ref{main-cor}.
\\

In Section~3, we study the Castelnuovo-Mumford regularity
of projective varieties as an application of $\S 2$.
In recent years upper bounds on the Castelnuovo-Mumford regularity
of a projective variety $X \subset \PN$ have been given by several
authors in terms of $\dim (X)$, $\deg (X)$, $\codim (X)$ and $k(X)$,
see, e.g., \cite{HM,MV,NS2}, where $k(X)$ is the 
Ellia-Migliore-Mir\'o Roig number measuring the deficiency
module, or sometimes called as the Rao module, see $\S 3$ for 
the definition.
A regularity bound
$\reg (X)  \le  \lceil (\deg (X) - 1)/\codim (X) \rceil
+ \max \{ k(X) \dim (X), 1 \}$
is known for a nondegenerate projective variety $X$,
see \cite{MV,NS2}.
Conversely, under the assumption that
a nondegenerate projective variety $X$ is ACM, that is,
the coordinate ring of $X$ is Cohen-Macaulay, if
$\reg (X)  \le  \lceil (\deg (X) - 1)/\codim (X) \rceil + 1$
and $\deg (X)$ is large enough,
then $X$ is a variety of minimal degree, see \cite{N1,TV}.
Moreover,  there gives a classification of
nondegenerate projective non-ACM varieties $X$ attaining
a regularity bound
$\reg (X)  =  \lceil (\deg (X) - 1)/\codim (X) \rceil
+ k(X) \dim (X)$.
In \cite{M},
under the assumption that $\deg (X)$ is large enough and
${\rm char}(K) =  0$, it is shown that
a projective non-ACM variety having the equality 
must be a curve on a rational
ruled surface, that is, on a Hirzeburch surface.
In $\S 3$, we show the corresponding result in the
positive characteristic case as an application of 
(\ref{main-th}), see Theorem~\ref{sharp-bound-th}.
\\

The first named author was partially supported
by MURST and GNSAGA of CNR (Italy).
The second named author was partially supported by
 Grant-in-Aid
for Scientific Research (no. 09740042),
Ministry of Education, Science, Sports and
Culture (Japan).
\\

\section{Regularity of a Generic Hyperplane \\
Section of Projective Curves in Positive \\
Characteristic}

Let $K$ be an algebraically closed field with
${\rm char}(K)  =  p  >  0$.  \\

In this section we will show that if $S  \subset \PN$
is a generic hyperplane section of an integral curve
with $\deg (S)$ large enough, then
either $S$ is in uniform position or 
$i(S)  \le  \lceil  (\deg (S)  -  1)/N  \rceil  -  1$.
Here the index of regularity
$i(S)$  of  $S$  is defined as the smallest integer $t$ such that
${\rm H}^1(\PN, {\cal I}_S(t))  =  0$.
(Notice that $\reg (S)  =   i(S) +  1  =  a(R) +  2$,
where $R$ is the coordinate ring of $S$ and $a(R)$
is an $a$-invariant of $R$, that is,
$a(R) = \max  \{  \ell \, |  \,  [ H_{\m_R}^1(R) ]_{\ell}
\ne  0  \}$ )  \\

A zero-dimensional scheme $S \subset \PN$ is called
in uniform position if ${\rm H}_Z(t) =
\max \{ \deg (Z), {\rm H}_S(t) \}$ for all $t$,
for any subscheme $Z$ of $S$, where ${\rm H}_Z$
and ${\rm H}_S$ denote the Hilbert function of $Z$
and $S$ respectively.
\\

A zero-dimensional scheme $S$, spanning $\PN$,
is called in linear semi-uniform position if there are
integers $v(i,S)$, simply written as $v(i)$,
$0 \le i  \le  N$ such that every
$i$-plane $L$ in $\PN$ spanned by linearly independent
$i+1$ points of $S$
contains exactly $v(i)$ points of $S$. 
A generic hyperplane section of a nondegenerate projective
integral curve is in linear semi-uniform position,
see \cite{B1}. 
We say $S$ is in linear general position if $v(i) = i +  1$
for all $i \ge 1$.
\\

Let $S$ be a zero-dimensional scheme of $\PN$ in linear
semi-uniform position. Then 
$v(i+1) \ge (v(1) - 1) v(i) + 1$ for $0 \le i \le N - 1$, see \cite{BC}.
Also, we have, by \cite{B1} (or see \cite{NS2}), $i(S)
\le  \lceil  (\deg (S)  -  1) / N  \rceil $.
\\

Further, we note that ``uniform position'' implies ``linear
general position'' and that ``linear general
position'' implies ``linear semi-uniform position''.
\\ 


\begin{lem}
\label{lem-v(1)=3}
Let $S \subset \PN$ be a generic hyperplane section
of a nondegenerate projective integral curve $C \subset {\PP}_K^{N+1}$
with $d = \deg (C)$.
Assume that  $N \ge 3$ $d \ge 25$.
If $v(1)  \ge 3$, then
$i(S)  \le  \lceil (d - 1)/N \rceil - 1$.
\end{lem}

\noindent {\it Proof}.
The assumption $v(1) \ge  3$ yields
$v(i) \ge  2^{i+1} - 1$ for  $0 \le  i \le  N$.
Put $v = v(N-1)$ and $w = v(N-2)$.
Note that $w \ge 2^{N-1} - 1$,
$v \ge (v(1) - 1) v(N-2) + 1 \ge  2 w + 1$
and $d \ge 2 v + 1 \ge 2^{N+1} - 1$.
\par

We have only to show that
${\rm H}^0({\cal O}_{\PN}(\ell))  \to
{\rm H}^0({\cal O}_S(\ell))$
is surjective, where $\ell =  \lceil (d - 1)/N \rceil  -  1$.
For any fixed point $P  \in  S$, we will show that there
is a union of $\ell$
hyperplanes $F = H(1)  \cup  \cdots \cup H(\ell)$ in $\PN$
such that
$S  \cap  F   =  S  \backslash \{ P \}$, as in the classical
Castelnuovo's method for finite sets in linear general position.
\\

First, let us take a hyperplane $H(1)$ which contains exactly
$v$ points of $S \backslash \{ P \}$ from the linear
semi-uniform position property. Then $H(1)$
does not contain $P$.
\par

Next, let us fix an $(N-2)$-plane $L$ in $H(1)$ such that
$L$ contains exactly $w$ points of $S \cap H(1)$.
Put $\ell_1 = \lfloor (d-v-1)/(v-w) \rfloor + 1$.
Now we will inductively construct hyperplanes $H(2),
\cdots, H(\ell_1)$ such that the number of points of
$(S \backslash \{ P \}) \cap (H(1) \cup \cdots \cup H(i))$
is $v + (i - 1) (v - w)$ for $i = 1, \cdots, \ell_1$.
In fact, since $d - 1 - v - (i - 1) (v - w) \ge v - w$
for $i \le \ell_1 - 1$, there exists a point $Q$ in $S 
\backslash (\{ P \} \cup
H(1) \cup \cdots \cup H(i))$ such that a hyperplane $M$ spanned
by $L$ and $Q$ does not contain $P$. 
Then $M$ contains exactly $v-w$ points of
$S \backslash (\{ P \} \cup H(1) \cup \cdots \cup H(i))$
from the linear semi-uniform position property.
So we take
$H(i+1) = M$.
Thus the union of $\ell_1$ hyperplanes $H(1) \cup \cdots
\cup H(\ell_1)$ contains $v + (\ell_1 - 1)(v - w)$ points
of $S$ and does not contain $P$.
Also, we note that
 $S \backslash (\{ P \} \cup H(1) \cup \cdots \cup H(\ell_1))$
consists of at most $v - w - 1$ points.
\par

However, we see that
$S \backslash (\{ P \} \cup H(1) \cup \cdots \cup H(\ell_1))$
consists of exactly $v - w - 1$ points.
In fact, if the number of the remaining points were
less than $v - w - 1$, then the hyperplane spanned by $M$
and a point from $S \backslash (H(1) \cup \cdots \cup H(\ell_1))$
would contain
at most $v - 1$ points of $S$, which contradicts with $v(N-1) = v$.
Thus we also have that there exist a hyperplane $G$
containing the $(N-2)$-plane $L$,
all the remaining points of
$S \backslash (\{ P \} \cup H(1) \cup \cdots \cup H(\ell_1))$
and the point $P$.
Of course $S \cap G$ consists of exactly $v$ points including
$P$.
\par

Since $S \cap G$ is in linear semi-uniform position in $G \cong 
\PP_K^{N-1}$, there are $\ell_2$ hyperplanes $M(\ell_1 + 1),
\cdots, M(\ell_2)$ of $\PP_K^{N-1}$ such that the union of them
contains the remaining points and does not contain $P$,
where $\ell_2 = \lceil (v - 1) / (N - 1) \rceil
(= \lfloor (v - 2)/(N-1) \rfloor + 1)$.
Thus we can take $\ell_2$ hyperplanes $H(\ell_1 + 1),
\cdots, H(\ell_2)$ of $\PN$ as desired.
Note that we used a fact from \cite{B1} 
that  ${\rm H}^0({\cal O}_{\PP_K^{N-1}}(t))  \to
{\rm H}^0({\cal O}_{S \cap G}(t))$
is surjective for all $t  \ge  \lceil (v - 1)/(N - 1) \rceil $, 
not necessarily for $t =  \lceil (v - 1)/(N - 1) \rceil - 1$,
without using the hypothesis of the induction on $N$.
So, if necessary, we may need to take a (possibly reducible)
hypersurface $F(1)$ of degree $\ell_2$ in place of the union of $\ell_2$
hyperplanes, and then go on the similar proof.
\par

Therefore we have $S \cap (H(1) \cup \cdots \cup H(\ell_1)
\cup \cdots \cup H(\ell_1 + \ell_2)) = S \backslash \{ P \}$
(or $S \cap (H(1) \cup \cdots \cup H(\ell_1)
\cup F(1)) = S \backslash \{ P \}$).
\\

Thus the proof is reduced to an arithmetic question.
In other words, we need to prove 
$\ell_1 + \ell_2 \le \ell$, namely,
\[ \left\lceil \frac{d-1}{N} \right\rceil - 
\left\lfloor  \frac{d - v - 1}{v -w} \right\rfloor 
- \left\lfloor  \frac{v - 2}{N-1} \right\rfloor  
\ge 3.  \]
Moreover, from the above argument, we remark that
$d = v + \ell_1 (v - w)$.
\\

First, assume that $N \ge 5$. Since
$v - w \ge w + 1 \ge 4 (N-1)$, it suffices to show that
$(d-1)/N - (d-v-1)/4(N-1) - (v - 2)/(N - 1) \ge 3$.
In fact, we easily have this inequality
by reducing it to the case $d = 2 v + 1$.
Hence we proved the case $N \ge 5$.
\par

Second, assume that $N = 4$. The inequality
$\lceil (d  -  1)/4 \rceil
 - \lfloor (d  -  v  -  1)/(v - w) \rfloor
 - \lfloor (v - 2)/3 \rfloor  \ge 3$ holds except for the case
$(d , v, w) = (32, 15, 7)$ or $(33, 15, 7)$.
But both cases contradict with $d = v + \ell_1 (v - w)$.
Hence we proved the case $N = 4$.
\par

Finally, assume that $N = 3$.
Then we have
$\lceil (d - 1)/3 \rceil 
- \lfloor (d - v - 1)/(v - w) \rfloor
- \lfloor (v - 2)/2 \rfloor
\ge 3$ except for the case $w = 3$ and
$(d, v) = (25, 7), (25, 8), (25, 10), (25, 12),
(28, 7)$ under the condition $d \ge 25$.
But all the exceptional cases contradict with
$d = v + \ell_1 (v - w)$.
Hence we proved the case $N = 3$.

\hfill \qed \bigskip

\begin{lem}
\label{lem-v(2)=4}
Let $S \subset \PN$ be a generic hyperplane section
of a nondegenerate projective integral curve 
$C \subset {\PP}_K^{N+1}$ with $d  =  \deg (C)$.
Assume that $N  \ge  3$ and $d \ge 23$.
If $v(1)  =  2$ and $v(2) \ge 4$, then
$i(S)  \le  \lceil (d - 1)/N \rceil - 1$.
\end{lem}

\noindent {\it Proof}.
In fact, by \cite{B2}, the assumption in (\ref{lem-v(2)=4})
yields that $\deg (C)  =   2^k$ for some $k  \ge  N$
and $v(i,S) = 2^i$ for all $i \le  N-1$ since $d \ge 23$.
In particular, $v(N-1) = 2^{N-1}$ and $v(N-2) = 2^{N-2}$.
\\

First assume that $N \ge 5$.
Just by copying the proof of (\ref{lem-v(1)=3})
as in the Castelnuovo's method, we see that the proof is
reduced to show an inequality
$ \lceil (2^k - 1)/N \rceil
- \lfloor (2^k - 2^{N-1} -1)/(2^{N-1}-2^{N-2})  \rfloor
- \lfloor (2^{N-1} - 2)/(N - 1) \rfloor
\ge 3$, namely,
\[ \left\lceil \frac{2^k - 1}{N} \right\rceil \ge
 2^{k-N+2} - 1 + 
 \left\lceil \frac{2^{N-1} - 1}{N - 1} \right\rceil,
 \]
which is easily shown.
Hence we proved the case $N \ge 5$.
\\

Next assume that $N = 3$.
As in the classical Castelnuovo's method, we will take
a union of hyperplanes with containing $S$
and without containing $P$.
\par

First let us take a hyperplane $H(1)$ with containing
exactly 4 points of $S \backslash \{ P \}$.
\par

Now we will inductively construct hyperplanes $H(2),
\cdots, H(\ell_1)$
such that the number of points of
$(S \backslash \{ P \}) \cap G(i)$
is $4i$ for $i = 1, \cdots, \ell_1$,
where $\ell_1   =    2^{k-3}$ and
$G(i) = H(1) \cup \cdots \cup H(i)$.
For any $i = 1, \cdots, \ell_1 - 1$,
we will show that there
exists a hyperplane $H(i+1)$ with containing
exactly 4 points of $S \backslash (\{ P \}
\cup G(i))$.
In fact, take 2 points $Q_1$ and $Q_2$ in
$S \backslash (\{ P \} \cup G(i))$.
Then there exists a point $Q_3$ 
from $S \backslash (\{ P, Q_1, Q_2 \} \cup G(i))$
such that the hyperplane
spanned by $Q_1$, $Q_2$ and $Q_3$ does not contain
any points of
$S \cap (\{ P \} \cup G(i))$, since
the number of points of $S \backslash (\{ P, Q_1, Q_2 \}
\cup G(i))$ is larger than
that of $S \cap (\{ P \} \cup G(i))$.
\par

So the number of the remaining point of $S \backslash (\{ P \}
\cup H(1) \cup \cdots \cup H(\ell_1))$
is $2^{k-1} - 1$.
Next we will inductively construct hyperplanes $H(\ell_1 + 1),
\cdots, H(\ell_1 + \ell_2)$ for some $\ell_2
\le \lceil (2^{k-1})/3 \rceil $,
satisfying that
$S \backslash \{ P \} = S \cap 
(H(1) \cup \cdots \cup 
H(\ell_1 + \ell_2))$.
In fact, assume that we already take hyperplanes
$H(1), \cdots, H(i)$
for $i \ge \ell_1$ satisfying some suitable condition.
If the number of the remaining points of
$S \backslash (\{ P \} \cup G(i))$, where
$G(i) = H(1) \cup \cdots \cup H(i)$, is larger than 3,
we can take the hyperplane $H(i+1)$ spanned by appropriate 3 points
from $S \backslash (\{ P \} \cup G(i))$ such that
$H(i+1)$ does not contain $P$.
So the number of the points of
$S \cap (H(i+1) \backslash G(i))$ is at least 3, and
possibly 4.
If the number of the remaining points of
$S \backslash (\{ P \} \cup G(i))$ is 3, then
we take hyperplanes $H(i+1)$ and $H(i+2)$ such that
$H(i+1) \cup H(i+2)$ contains the remaining 3 points
of $S \backslash (\{ P \} \cup G(i))$ and does not contain
$P$. 
If the number of the remaining points of
$S \backslash (\{ P \} \cup G(i))$ is either 1 or 2, then
we take a hyperplane $H(i+1)$ such that
$H(i+1)$ contains the remaining 1 or 2 points
of $S \backslash (\{ P \} \cup G(i))$ and does not contain
$P$. 
\\

Thus the proof is reduced to an arithmetic question
as in (\ref{lem-v(1)=3}). Namely, $\ell_1 + \ell_2
\le  \lceil (2^k - 1)/3 \rceil - 1$, in other words,
\[ \left\lceil \frac{2^k - 1}{3} \right\rceil -
2^{k-3} - \left\lceil \frac{2^{k-1}}{3} \right\rceil 
\ge 1.
\]
Then we easily see the inequality except for the case 
$k = 3, 4$. 
\par

Hence we proved the case $N = 3$. 
\\

Finally assume that $N = 4$.
Again we will prove as in the classical Castelnuovo's method.
\par

First let us
take hyperplane $H(1)$ with containing exactly
8 points of $S \backslash \{ P \}$.
\par

Now we will inductively construct hyperplanes $H(2),
\cdots, H(\ell_1)$ for some integer $\ell_1
\le \lfloor (2^{k-1} + 1)/7 \rfloor$
such that $S \cap (H(i+1) \backslash G(i))$
contains at least 7 points and does not contain $P$,
where $G(i) = H(1) \cup \cdots \cup H(i)$.
In fact, take 2 points $Q_1$ and $Q_2$ from
$S \backslash (\{ P \} \cup G(i))$.
Then there exists a point $Q_3$ 
in $S \backslash (\{ P, Q_1, Q_2 \} \cup G(i))$
such that the 2-plane $L$
spanned by $Q_1$, $Q_2$ and $Q_3$ does not contain
any points of
$S \cap (\{ P \} \cup G(i))$ 
if the number of points of $S \backslash (\{ P, Q_1, Q_2 \}
\cup G(i))$ is larger than
that of $S \cap (\{ P \} \cup G(i))$.
In other words,
we can take such $L$ if $S \backslash 
(\{ P \} \cup G(i))$
contains at least $2^{k-1} + 2$ points.
Thus the 2-plane $L$ contains exactly 4 points
of $S \backslash (\{ P \}
\cup G(i))$, and we put $S \cap L = \{Q_1, \cdots, Q_4 \}$.
Then there exists a point $Q_5$
from $S \backslash (\{ P, Q_1, \cdots, Q_4 \} \cup G(i))$
such that the hyperplane $M$ spanned by the point $Q_5$
and the 2-plane $L$ 
contains at least two points of
$S \backslash (\{ P, Q_1, \cdots, Q_4 \} \cup G(i))$
without containing $P$,
if the number of points of $S \backslash (\{ P, Q_1, 
\cdots, Q_4 \} \cup G(i))$ minus 2 is larger than
that of $S \cap (\{ P \} \cup G(i))$.
In this case we put $H(i+1) = M$.
In other words,
we can go on this process if $S \backslash 
(\{ P \} \cup G(i))$
contains at least $2^{k-1} + 4$ points.
Thus we constructed a union of hyperplanes
$G(\ell_1) = H(1) \cup \cdots \cup H(\ell_1)$
such that $G(\ell_1)$ contains at least
$2^{k-1} - 4$ points of $S$ and does not contain
$P$ for some $\ell_1  \le \lfloor (2^{k-1} + 1)/7 \rfloor $.
\par

So the number of the remaining point of $S \backslash (\{ P \}
\cup H(1) \cup \cdots \cup H(\ell_1))$
is at most $2^{k-1} + 3$.
Next we will inductively construct hyperplanes $H(\ell_1 + 1),
\cdots, H(\ell_1 + \ell_2)$ for some integer $\ell_2  \le 
 2^{k-3}+2 $
satisfying that
$S \backslash \{ P \} = S \cap 
(H(1) \cup \cdots \cup 
H(\ell_1 + \ell_2))$.
Assume that we already take hyperplanes $H(1), \cdots, H(i)$
for $i \ge \ell_1$ satisfying some suitable condition.
If the number of the remaining points of
$S \backslash (\{ P \} \cup G(i))$, where
$G(i) = H(1) \cup \cdots \cup H(i)$, is larger than 6,
we can take a hyperplane $H(i+1)$ with containing at least 4 points
of $S \backslash G(i)$ and without containing $P$.
So the number of $S \cap (H(i+1) \backslash G(i))$ is at least 4,
and possibly more.
If the number of the remaining points of
$S \backslash (\{ P \} \cup G(i))$ is 6, then
we take hyperplanes $H(i+1)$, $H(i+2)$ and $H(i+3)$ with
$H(i+1) \cup H(i+2) \cup H(i+3)$ containing the
 remaining 6 points
of $S \backslash (\{ P \} \cup G(i))$ and without containing
$P$.
If the number of the remaining points of
$S \backslash (\{ P \} \cup G(i))$ is either 3, 4 or 5, then
we take hyperplanes $H(i+1)$ and $H(i+2)$ with
$H(i+1) \cup H(i+2)$ containing the remaining 3, 4 or 5 points
of $S \backslash (\{ P \} \cup G(i))$ and without containing
$P$. 
If the number of the remaining points of
$S \backslash (\{ P \} \cup G(i))$ is either 1 or 2, then
we take a hyperplane $H(i+1)$ with
containing the remaining 1 or 2 points
of $S \backslash (\{ P \} \cup G(i))$ and without containing
$P$. 
Thus we see that there exist hyperplanes $H(\ell_1 + 1),
\cdots, H(\ell_1 + \ell_2)$ as desired.\\

Thus the proof is reduced to an arithmetic question
as in (\ref{lem-v(1)=3}). Namely,
$\ell_1 + \ell_2 \le \lceil (2^k - 1)/4 \rceil - 1$,
in other words,
\[ \left\lceil \frac{2^k - 1}{4} \right\rceil -
\left\lfloor \frac{2^{k-1} +  1}{7} \right\rfloor -
2^{k-3}
\ge 1.  \]
Then we easily see the inequality.
\par

Hence we proved the case $N = 4$. 
\hfill \qed \bigskip

\begin{th}
\label{main-th}
Let $S \subset \PN$ be a generic hyperplane section of
a nondegenerate projective integral
curve $C \subset {\PP}_K^{N+1}$  with $d = \deg (C)$.
If
$d  \ge  \max \{ N^2  +  2 N  +  2,  25  \}$
and
$i(S)    =    
\lceil  (d   -   1)/N   \rceil$,
then $S$ lies on a rational normal curve.
\end{th}

\noindent {\it Proof}.
For the case $N  =  2$, the corresponding
result as in \cite[(3.2)]{Y1} on the $h$-vector for the positive
characteristic case is true, see \cite[(1.1)]{BY} or
\cite{GM,HTV}.
So the assertion follows from the proof of \cite[(2.5)]{M}.
\par

We may assume that 
$N  \ge  3$ and that the Uniform Position Lemma fails for
the curve $C$. Note that $d \ge 25$. Then, by \cite[(2.5)]{R},
 $C$ satisfies either
(i) every secant of $C$ is a multisecant, that is,
$v(1) \ge 3$, or (ii) every plane spanned by three points
contains one more point of $C$, that is,
$v(1) = 2$ and $v(2) \ge 4$.
Therefore, by (\ref{lem-v(1)=3}) and (\ref{lem-v(2)=4}),
we obtain that $i(S)    \le    
\lceil  (d   -   1)/N   \rceil - 1$.
So we exclude the case.
\par

Hence the assertion is proved.
\hfill \qed \bigskip

\begin{lem}
\label{prop-N=2}
Let $S \subset \PP_K^2$ be a generic hyperplane section
of a nondegenerate integral space curve $C$ with $d = \deg (C)$.
If $v(1) \ge 4$, then
$i(S)  \le  \lceil (d - 1)/2  \rceil  - 1$.
\end{lem}

\noindent {\it Proof}.
Put $v = v(1)$.
Following the Castelnuovo's method, we will have the corresponding
proof as in (\ref{lem-v(1)=3}).
For any fixed point $P  \in  S$, we have only to show that there
is a union of $\ell$
lines $F = L(1)  \cup  \cdots \cup L(\ell)$ in $\PN$
such that
$S  \cap  F   =  S  \backslash \{ P \}$, where
$\ell = \lceil (d - 1)/2 \rceil - 1$.
\\

First, let us take a line $L(1)$ which contains exactly
$v$ points of $S \backslash \{ P \}$ from the linear
semi-uniform position property. Then $L(1)$
does not contain $P$.
\par

Next, let us fix a point $Q$ of $L(1)$ and
put $\ell_1 = \lfloor (d - v - 1)/(v - 1) \rfloor$.
Then we can construct lines $L(2),
\cdots, L(\ell_1)$, by taking inductively
a line $L(i+1)$ with containing $Q$
and without containing any points of $(\{ P \} \cup
L(1) \cup \cdots \cup L(i)) \backslash \{ Q \}$
for $1 \le i \le \ell_1-1$.
\par

Moreover, since $S \backslash (\{ P \} \cup L(1) \cup \cdots \cup L(\ell_1))$
consists of at most $v - 2$ points (and in fact exactly $v -2$ points),
 we can take appropriate
$v-2$ lines $L(\ell_1+1), \cdots, L(\ell_1 + v - 2)$
with containing the remaining points of $S \backslash \{ P \}$
and without containing $P$. 
\\

Thus the proof is reduced to an arithmetic question.
In other words, 
$\ell_1 + v - 2 \le \ell$, namely,
$\lceil (d-1)/2 \rceil - \lfloor 
(d - v - 1)/(v -1)  \rfloor - v + 1 
\ge 0$,
which is easily shown.
\par

Hence the assertion is proved.
\hfill \qed \bigskip

\begin{lem}
\label{prop2-N=2}
Let $S \subset \PP_K^2$ be a generic hyperplane section
of a nondegenerate integral space curve $C$ with $d = \deg (C)$.
If $v(1) = 3$ and $d \ge 24$, then
$i(S)  \le  \lceil (d - 1)/2  \rceil  - 1$.
\end{lem}

\noindent {\it Proof}.
Following the Castelnuovo's method, we will have the corresponding
proof as in (\ref{lem-v(2)=4}), the case $N=3$.
For any fixed point $P  \in  S$, we have only to show that there
is a union of $\ell$
lines $F = L(1)  \cup  \cdots \cup L(\ell)$ in $\PN$
such that
$S  \cap  F   =  S  \backslash \{ P \}$, where
$\ell = \lceil (d - 1)/2 \rceil - 1$.
\\

First, let us take a line $L(1)$ which contains exactly
3 points of $S \backslash \{ P \}$ from the linear
semi-uniform position property. Then $L(1)$
does not contain $P$.
\par

Put $\ell_1 = \lfloor (d - 4)/6 \rfloor + 1$.
Then we can construct lines $L(2),
\cdots, L(\ell_1)$, by taking inductively
a line $L(i+1)$ without containing any points of 
$\{ P \} \cup
L(1) \cup \cdots \cup L(i)$
for $1 \le i \le \ell_1-1$.
\par

Moreover, since $S \backslash (\{  P \} \cup L(1) \cup \cdots \cup L(\ell_1))$
consists of at most $\lceil (d + 1)/2 \rceil$ points, we can take
appropriate $\ell_2$
lines $L(\ell_1 + 1), \cdots, L(\ell_2)$
with containing the remaining points of 
$S \backslash \{ P \}$ and without containing $P$,
where $\ell_2 = \lceil (d + 3) / 4 \rceil$.
\\

Thus the proof is reduced to an arithmetic question.
In other words, 
$\ell_1 + \ell_2 \le \ell$, namely,
$\lceil (d-1)/2 \rceil - \lfloor 
(d - 4)/6 \rfloor - \lceil (d+3)/4 \rceil 
\ge 2$, 
which is easily shown for $d \ge 24$.
\par

Hence the assertion is proved.
\hfill \qed \bigskip

We say that a nondegenerate projective integral curve $C$
is very strange if a generic hyperplane section $S$ of $C$
is not in linear general position.

\begin{prop}
\label{main-cor}
Let $S \subset \PN$ be a generic hyperplane section
of a nondegenerate projective integral curve $C \subset {\PP}_K^{N+1}$.
Assume that $C$ is very strange.
If $d  =  \deg (C)  \ge 25$, then
$i(S)  \le  \lceil  (d - 1)/N  \rceil  - 1$ \\
\end{prop}

\noindent {\it Proof}.
It immediately follows from (\ref{lem-v(1)=3}), (\ref{lem-v(2)=4}),
 (\ref{prop-N=2}), (\ref{prop2-N=2}) and the proof of (\ref{main-th}).
\hfill \qed \bigskip

\section{An Application to a Sharp Bound on the \\
Castelnuovo-Mumford Regularity}

Let $K$ be an algebraically closed field.
Let $S  =  K[x_0,  \cdots,  x_N]$ be the polynomial ring
and $\m  =  (x_0,  \cdots,  x_N)$ be the irrelevant ideal.
Let $X$ be a projective scheme of $\PN = \Proj (S)$.
For an integer $m$,
$X$ is said to be $m$-regular if 
${\rm H}^i({\PP}_K^N, {\cal I}_X(m-i))  =  0$
for all $i  \ge  1$.
The Castelnuovo-Mumford regularity
of $X  \subset   \PN$ is the least such $m$
and is denoted by $\reg (X)$. \\

Let $k$ be a nonnegative integer.
Then $X$ is called $k$-Buchsbaum if
the graded $S$-module
${\rm M}^i(X)
  =
  \oplus_{\ell \in {\Bbb Z}}{\rm H}^i({\Bbb P}_K^N, {\cal I}_X(\ell))$,
called the deficiency module of $X$,
is annihilated by $\m^k$
for  $1    \le      i   \le     \dim (X)$,
see, e.g., \cite{Mi,MM}. 
On the other hand, $X$ is called strongly $k$-Buchsbaum
if $X \cap V$ has the $k$-Buchsbaum property for any
complete intersection $V$ of $\PN$ with $\codim (X \cap V)
= \codim (X) + \codim (V)$, possibly $V = \PN$.
So ``strongly $k$-Buchsbaum'' implies ``$k$-Buchsbaum''.
Further we call the minimal nonnegative integer $n$, if there exists,
such that $X$ is $n$-Buchsbaum (resp. strongly $n$-Buchsbaum),
as the Ellia-Migliore-Mir\'o Roig number (resp. the strongly
Ellia-Migliore-Mir\'o Roig number)
of $X$ and denote by $k(X)$ (resp. $\bar{k}(X)$), see \cite{M}.
In case $X$ is not $k$-Buchsbaum for all $k  \ge   0$,
then we put $k(X)   =   \bar{k}(X) =   \infty$.
Note that $k(X)  <  \infty$ if and only if $\bar{k}(X) < \infty$.
Moreover it is equivalent to saying that 
$X$ is locally Cohen-Macaulay
and equi-dimensional.
\\

Upper bounds on the Castelnuovo-Mumford regularity
of a projective variety $X$ are given in terms of
$\dim (X)$, $\deg (X)$, $\codim (X)$, $k(X)$ and $\bar{k}(X)$.
Moreover, in case ${\rm char}(K)  =  0$, the extremal
cases for the bounds are classified under a certain
assumption.

\begin{prop}
\label{main-prop}
Let $X$ be a nondegenerate projective
variety in ${\PP}_K^N$. 
Assume that $X$ is not ACM, that is, $k(X) \ge 1$.
Then 
\medskip

\noindent
$(a) \quad \reg (X) \le \lceil
(\deg (X) - 1)/\codim (X)  \rceil + k(X) \dim(X)$.

\medskip

\noindent
$(b) \quad \reg (X) \le \lceil
(\deg (X) - 1)/\codim (X)  \rceil + \bar{k}(X) \dim(X) -
\dim (X) + 1$.

\medskip

Furthermore, assume that ${\rm char}(K) = 0$ and
$\deg (X)  \ge  2 \, \codim (X)^2  +  \codim (X)  +  2$.
If the equality $\reg (X) = \lceil
(\deg (X) - 1)/\codim (X) \rceil + k(X) \dim(X)$ holds,
then $X$ is a curve on a rational ruled surface.
\end{prop}

\noindent {\it Proof}.
See \cite{M,MV,NS2}.
\hfill \qed \bigskip

Now we will study the extremal case for the
inequality in (\ref{main-prop}) in positive characteristic
case. We assume that a variety is not ACM,
see \cite{N1} for the ACM case.

\begin{th}
\label{sharp-bound-th}
Let $X$ be a nondegenerate projective
variety in ${\PP}_K^N$ with $k(X) \ge 1$.
Assume that either ${\rm char}(K) = 0$ and
$\deg (X)  \ge   \codim (X)^2  +  2 \, \codim (X)  +  2$,
or ${\rm char}(K)  = p > 0$ and 
$\deg (X)  \ge  \max \{
2 \, \codim (X)^2  +  \codim (X)  +  2, 25 \}$.

\medskip

\noindent
$(a)$
If the equality $\reg (X) = \lceil
(\deg (X) - 1)/\codim (X) \rceil + k(X) \dim(X)$ holds,
then $X$ is a curve on a rational ruled surface.

\medskip

\noindent
$(b)$
If the equality $\reg (X) = \lceil
(\deg (X) - 1)/\codim (X) \rceil + \bar{k}(X) \dim(X)
-\dim (X) + 1$ holds,
then $X$ is a curve on a rational ruled surface.

\end{th}

\noindent {\it Proof}.
We will prove (a). The proof of (b) is similar as in (a), which is
left to the readers. \\

First we assume that 
${\rm char}(K)  = p > 0$ and 
$\deg (X)  \ge  \max \{
2 \, \codim (X)^2  +  \codim (X)  +  2, 25 \}$.
The lemmas (2.5), (2.6), (2.7) and (2.8) in \cite{M}
work for the case ${\rm char}(K) = p > 0$, although
an assumption ${\rm char}(K) = 0$ is mentioned in \cite{M}.
However, for the positive characteristic case, we cannot
apply \cite[(2.5)]{M} as an inductive step, because
a generic hyperplane section of an integral curve is not
necessarily in uniform position.
In other words, the corresponding proof as in \cite{M} 
works for the
positive characteristic case, except for
the Uniform Position Lemma.
\par

Thus, by applying Theorem~\ref{main-th} in place of \cite[(2.5)]{M},
we have the assertion.
\\
 
On the other hand, for the case ${\rm char}(K) = 0$ and
$\deg (X)  \ge   \codim (X)^2  +  2 \, \codim (X)  +  2$,
we use \cite[(3.3)]{N2} in place of \cite[(2.6),(2.8)]{M}.
(Notice that \cite[(3.3)]{N2} is a consequence of 
the ``Socle Lemma'',
see \cite{HU}, and cannot be applied for the positive
characteristic case.)
Hence we have the assertion.
\hfill \qed 

\end{document}